\documentclass{amsart}

\usepackage{amssymb}

\newtheorem{thm}{Theorem}[section]

\newtheorem{lem}[thm]{Lemma}
\newtheorem{cor}[thm]{Corollary}
\theoremstyle{definition}

\theoremstyle{remark}
\newtheorem{remark}[thm]{Remark}
\numberwithin{equation}{section}

\begin{document}

\title[Higher Order]{On Higher Order Elliptic Differential Operators I: The Geometric Problem}

\author{David Raske}

\address{David Raske, 3635 Washtenaw Ave, Ann Arbor, Michigan, 48104, United States}

\email{nonlinear.problem.solver@gmail.com}



\begin{abstract} In this note we will examine the issue of what types of phenomena can be associated with higher order elliptic differential operators on manifolds but not with second order elliptic differential on manifolds. Then we will discuss the similarities between coercive higher order elliptic differential operators on manifolds and coercive second order elliptic differential operators on manifolds.
\end{abstract}

\maketitle
\section{Introduction}

The main body of this paper is divided into two sections. The first section is devoted to the differences between higher order elliptic differential operators on manifolds and second order differential operators on manifolds. In particular, we will see that there are a lot of elliptic differential operators on manifolds such that the lowest eigenvalue is not simple and such that a Green's operator exists and is not positivity preserving. The second section is devoted to the similarities between higher order elliptic differential operators on manifolds and second order elliptic differential operators on manifolds. In particular, we will see that if $(M,g)$ is a closed, smooth Riemannian manifold and G is a Green's function associated with a coercive elliptic differential operator on $(M,g)$, then $\int \int_{M \times M} G(x,y) dx dy > 0$, and the Green's function is positive on the diagonal. In short, higher order elliptic differential operators on manifolds can be very different than second order elliptic operators on manifold, but nonetheless there will always be similarities when it comes to properties of the Green's functions associated with coercive elliptic differential operators on closed, smooth manifolds. 

\section{Higher order phenomena}
\subsection{Some higher order elliptic differential operators on smooth manifolds have a non-simple lowest eigenvalue.}

The lowest eigenvalue of a formally self-adjoint second order elliptic operator on a smooth, closed manifold $M$ has to be simple, but the lowest eigenvalue of a higher order formally self-adjoint elliptic operator on a smooth, closed manifold does not have to be simple. In particular we have

\begin{thm}
Let $(M,g)$ be a smooth, closed Riemannian manifold of dimension greater than or equal to one. Let $m$ be a positive integer. Let $P_g:  C^\infty(M) \rightarrow C^\infty(M)$ be a $2m$th order formally self-adjoint elliptic differential operator on $M$ that has a lowest eigenvalue. Let $\Sigma$, the real spectrum of $P_g$, be $\{\lambda_i\}_{i=1}^ \infty$, where $\lambda_1 \leq \lambda_2 \leq \lambda_3 \leq ...$, and where we repeat each eigenvalue according to its finite multiplicity. Let $k$ be a positive integer such that $\lambda_k < \lambda_{k+1}$. Then the self-adjoint elliptic differential operator $(P_g-\frac{\lambda_k + \lambda_{k+1}}{2})^2$ has the property that its lowest eigenvalue has a multiplicity of at least two.
\end{thm}

\begin{proof} First note that the spectral theorem gives us that there exists a sequence of smooth real valued functions $\{e_i\}_{i=1}^\infty$ on $M$ such that $P_g e_i = \lambda_i e_i$ for all positive integers $i$. Furthermore the spectral theorem allows us to assume without loss of generality that the sequence $\{e_i\}_{i=1}^\infty$ forms an orthonormal basis for $L^2(M)$, so let us do so.

Now let $k$ be a positive integer such that $\lambda_k < \lambda_{k+1}$ and note that we have
$$
(P_g-\frac{\lambda_k + \lambda_{k+1}}{2})^2 e_i = (\lambda_i - \frac{\lambda_k + \lambda_{k+1}}{2})^2 e_i
$$
for all positive integers $i$. Note as well that algebra gives us that \
$$
(\lambda_k - \frac{\lambda_k + \lambda_{k+1}}{2})^2 = (\lambda_{k+1} - \frac{\lambda_k + \lambda_{k+1}}{2})^2 \leq (\lambda_i - \frac{\lambda_k + \lambda_{k+1}}{2})^2
$$
for all positive integers $i$. 

Now suppose that there exists a real number $\gamma$ that is lower than $(\lambda_k - \frac{\lambda_k + \lambda_{k+1}}{2})^2$ and a not identically zero function $w \in C^\infty(M)$ such that $(P_g-\frac{\lambda_k + \lambda_{k+1}}{2})^2 w = \gamma w$. It would follow that

\begin{equation}
\begin{split}
\gamma \int_M w e_i dv_g & = \int_M (P_g-\frac{\lambda_k + \lambda_{k+1}}{2})^2 w e_i dv_g \\ &= \int_M w (P_g-\frac{\lambda_k + \lambda_{k+1}}{2})^2 e_i dv_g \\ &= (\lambda_i - \frac{\lambda_k + \lambda_{k+1}}{2})^2 \int_M w e_i dv_g.
\end{split}
\end{equation}
Note that this forces $w$ to be orthogonal, in the $L^2$ sense, to all the $e_i$, This produces a contradiction, since $\{e_i\}_{i=1}^\infty$ forms a complete basis for $L^2(M)$. It follows that the lowest eigenvalue of $(P_g- \frac{\lambda_k + \lambda_{k+1}}{2})^2$ is $(\lambda_k - \frac{\lambda_k + \lambda_{k+1}}{2})^2$ and we established  above that this eigenvalue has an eigenspace with a dimension greater than one. The theorem follows. 
\end{proof}

We also have

\begin{thm}
Let $M$ be a smooth, closed Riemannian manifold of dimension greater than or equal to one. Let $P: C^\infty(M) \rightarrow C^\infty(M)$ be a formally self adjoint elliptic differential operator. Suppose that $P$ has a lowest eigenvalue, and suppose that $P$ has an eigenvalue $\lambda$ that has a multiplicity greater than one. Then the lowest eigenvalue of $(P-\lambda)^2$ has a multiplicity greater than one.
\end{thm}

\begin{proof} Suppose that there exists a real number $\lambda$ and two smooth functions $u: M \rightarrow \mathbb{R}$ and $v: M \rightarrow \mathbb{R}$ such that $Pu = \lambda u$ and $Pv = \lambda v$. This in turn, allows us to write $(P- \lambda )u=0$ amd $(P-\lambda)v=0$.  We now have that $(P-\lambda)^2 u=0$ and $(P-\lambda)^2 v=0$, and hence $u$ and $v$ are in the kernel of $(P-\lambda)^2$. Since the eigenvalues of $(P-\lambda)^2$ are non-negative, we have the theorem.
\end{proof}

\subsection{More on signchanging eigenfunctions associated with the lowest eigenvalue of an elliptic differential operator}

The eigenfunction associated with the lowest eigenvalue of a second order formally self-adjoint elliptic differential operator on a smooth, closed manifold has to be signed. But an eigenfunction associated with the lowest eigenvalue of a higher order formally self-adjoint elliptic operator on a smooth, closed, Riemannian manifold does not have to be signed. In particular we have

\begin{thm}
Let $(M,g)$ be a closed, smooth Riemannian manifold of dimension one or higher. Let $P: C^\infty(M) \rightarrow C^\infty(M)$ be a formally self adjoint elliptic differential operator with the property that it has a lowest eigenvalue. Suppose that the lowest eigenvalue is not simple. Then $P$ has a a sign changing eigenfunction associated with its lowest eigenvalue.
\end{thm}

\begin{proof} Suppose all of the eigenfunctions associated with the lowest eigenvalue of $P$ are signed. It would follow that we can not find an orthogonal basis for the eigenspace of the lowest eigevalue. This produces a contradiction because the Gram-Schmidt procedure guarantees that we can find an orthogonal basis for the eigenspace for the lowest eigevalue of $P$.
\end{proof} 

Let $(M,g)$ be a smooth, closed Riemannian manifold of dimension one or higher. Let $P: C^\infty(M) \rightarrow C^\infty(M)$ be a formally self-adjoint elliptic differential operator with a lowest eigenvalue, We will say that an elliptic differential operator $P: C^\infty(M) \rightarrow C^\infty(M)$ is a type 1 FSA elliptic differential operator if the lowest eigenvalue is simple and the eigenfunctions associated to it are signed. We will say that a formally self-adjoint elliptic differential operator on a smooth manifold, $P$, is a type 2 FSA elliptic differential operator if it is not a type 1 FSA  elliptic differential operator. Note that whether or not an operator is a type 1 FSA operator is an important issue.  Now, let $\Delta_g: C^\infty(M) \rightarrow C^\infty(M)$ be defined as follows: $\Delta_g u := -\text{div}_g (\nabla u)$. If we restrict our attention to fourth order elliptic differential operators on manifolds we can find a large family of type 2 operators just knowing the first non-zero eigenvalue of $\Delta_g$. To see why this is the case we will now prove

\begin{thm} Let $(M,g)$ be a smooth, closed Riemannian manifold of dimension $m \geq 1$. Let $\lambda_2$ be the second lowest eigenvalue of $\Delta_g$.  Let $A$ be a smooth, symmetric $(0,2)$-tensor field on $M$. Let $P_g$ be the elliptic differential operator on $M$ that is given on functions by $P_g u = \Delta^2_g  u - \text{div}_g(A(\nabla u)^\sharp)$ . Then, $P_g$ has a sign-changing  eigenfunction that is associated with the lowest eigenvalue of $P_g$, if $A$ has the property that $\max_{x \in M} \gamma_m(x) \leq -\lambda_2$ where $\gamma_m(x)$ is 
$$
\max_{w \in T_x(M) \setminus \{0\}}  \frac{A(x)(w,w)}{g(x)(w,w)}.
$$
($A(x)$ is  $A$ at $x \in M$, and $g(x)$ is $g$ at $x \in M$.) This, in turn, allows us to conclude that if $A$ has the property that $\max_{x \in M} \gamma_m(x) \leq -\lambda_2$, then $P_g$ is a fourth order Type 2 FSA elliptic differential operator.
\end{thm}

A proof of this theorem can be found in [Raske 2026].

\subsection{There exists a coercive higher order elliptic differential operator on a smooth manifold that has a simple lowest eigenvalue that is associated with a signed eigenfunction and such that the operator’s Green’s function is not signed. }

 First we prove a result that stands out here because it seems like a result that a second order elliptic differential operator would have. We will see that, nonetheless, it will have higher order consequences.

\begin{lem}
Let $(M,g)$ be  closed, smooth Riemannian manifold of dimension one or higher. Let $k$ be a positive integer, and let $L$ be an $2k$th order, elliptic differential operator on $M$. Suppose as well that the kernel of $L$ consists of constant functions, and  $\int_M uLu \, dvol_g \geq 0$ for all $u \in H^{k}(M)$. Finally, let $\lambda$ be a real positive number. Then we have that the lowest eigenvalue of $L + \lambda$ is simple and positive. Furthermore non-zero constant functions are the eigenfunctions associated with the lowest eigenvalue of $L +\lambda$.
\end{lem}
 
\begin{proof} Let $F: H^k(M) \setminus \{0\}\rightarrow \mathbb{R}$ be the Rayleigh quotient
$$
\frac{\int_M (u L u + \lambda u^2) \, \text{dvol}_g}{\int_M u^2 \text{dvol}_g.}
$$
We will use the calculus of variations to show that $F$ is minimized over $H^k(M) \setminus \{0\}$ by non-zero constant functions. First note that since $F$ is homogeneous we can assume that $\int_M u^2 \, \text{dvol}_g =1$. We can then use our hypothesis that $\int_M w Lw \, \text{dvol}_g \geq 0$ for all $w \in H^k(M)$ to establish that $F$ is bounded below on $H^k(M) \setminus \{0\}$. We can now use our hypothesis that the kernel of $L = \{constants\}$ to establish that the minimizer of $F$ over $H^k(M) \setminus \{0\}$ exists and is achieved by constants. We know that $F$ is not minimized over $H^k(M) \setminus \{0\}$ by a nonconstant function $z$ because $\int zLz \, \text{dvol}_g > 0$ if $z$ is not a constant. The lemma follows.
\end{proof}

\begin{remark}
 Let $(M,g)$ be a smooth, closed Riemannian manifold. Let $\lambda$ be a positive real number. Let $P: C^\infty(M) \rightarrow C^\infty(M)$ be the elliptic differential operator $\Delta_g^2 + \lambda$. Then the above lemma allows us to conclude that the lowest eigenvalue of $L + \lambda$ is simple and positive. Furthermore non-zero constant functions are the eigenfunctions associated with the lowest eigenvalue of $L +\lambda$. This fact will be useful in the sequel.
\end{remark}

Now we need an ODE result. It was communicated to the author in [Yang 2009].

\begin{lem} Let $u: \mathbb{R} \rightarrow \mathbb{R}$ be defined as follows: $u(x)= (\sin{x})^4 + \epsilon$. Then for small negative choices of $\epsilon$, $u_{xxxx} + 944 u$ will be positive and periodic with wavelength of $2\pi$.
\end{lem}

\begin{proof} First note that $u_{xxxx} = 1200(\cos{x})^4-2160(\cos{x})^2 + 984 + 944 \epsilon$. Now, note that for small negative choices of $\epsilon$, $u_{xxxx} + 944 u \geq 0$. But,we also have that $u < 0$ somewhere, because $u=(\sin{x})^4 + \epsilon$. Finally note that $u_{xxxx} + 944u$ is periodic with wavelength $2\pi$.
\end{proof}

Finally, we can state and prove

\begin{thm} Let $(M,g)$ be the closed, smooth Riemannian manifold $(S^1,g_{can})$. Then the elliptic differential operator $\Delta_{g_{can}}^2 + 944$ will have a positive simple lowest eigenvalue. The eigenfunctions associated with the lowest eigenvalue will be nonzero constants. Finally, the Green's function of $\Delta_{g_{can}}^2 + 944$ will change sign.
\end{thm}

\begin{proof} Combine the above two lemmas and recall that a Green's function is sign-changing if the corresponding Green's operator is not positivity preserving.
\end{proof}

\section{Similarities between coercive second order elliptic differential operators and coercive higher order elliptic differential operators}
\subsection{A Green's functions associated with a positive Green's operators has to be positive on the diagonal}

Before we proceed,to prove results about Green's functions, we need to define some Hilbert spaces. Let $(M,g)$ be a smooth, closed Riemannian manifold of dimension one or greater. Let $k$ be a positive integer. Then the Hilbert space $H^{k}(M)$ is the set of $u \in L^2(M)$ such that $Pu =f \in L^2(M)$ (in the weak sense) whenever $P$ is a differential operator of order $\leq k$. $H^{-k}(M)$ is the continuous dual of $H^k(M)$. 

We can now prove results about Green's functions on $M \times M$. Let $M$ be a closed smooth Riemannian manifold of dimension $n$. Let $k$ be a positive integer, and let $P: H^{2k}(M) \rightarrow L^2(M)$ be a $2k$th order coercive classical elliptic differential operator. Let $L: H^k(M) \rightarrow L^2(M)$ be a positive $k$th order elliptic differential operator on a smooth manifold such that $L^2 = P$.  We will call $L = P^{1/2}$ and we will define $B: H^k(M) \times H^k(M) \rightarrow \mathbb{R}$ as follows: $B(u,v) = \int_M P^{1/2} u P^{1/2} v  dv_g$.  

\begin{lem} $B$ is an inner product for $H^k(M)$ and $(B(u,u))^{1/2}$ is a norm on $H^k(M)$ equivalent to the standard one.
\end{lem}

\begin{proof} The only thing that we need that is not obvious is that the norm $(B(u,u))^{1/2}$ is a norm on $H^k(M)$ that is equivalent to the standard one. This was proven in [Seeley 1967], so we have the lemma.
 \end{proof}

Now, fix $x \in M$ and let $\delta_x$ be the delta function supported on $x$. Suppose $2k > n$. Then we have 

\begin{lem} Let $(M,g)$ be a closed smooth Riemannian manifold of dimension $n$. Let $k$ be a positive integer such that $2k > n$. Let $x \in M$. Then $\delta_x \in H^{-k}(M)$.
\end{lem}

\begin{proof} First note that if $2k > n$, then functions in $H^{k}(M)$ are in $C(M)$, the space of continuous functions on $M$. Then, since $R(M)$, the space of Radon measures on $M$, is the continuous dual of $C(M)$, it follows that it is contained in $H^{-k}(M)$, the continuous dual of $H^k(M)$. Now, let us note $\delta_x$, the delta function supported at $x$, is a Radon measure. The lemma follows.
\end{proof}

We will now define a weak solution of $P u = f$, where $f \in H^{-k}(M)$. We will define it as a solution of $B(u,v) = <f,v>$ for all $v \in H^k(M)$, where $<\cdot,\cdot>$ is the duality paring between $H^{-k}(M)$ and $H^{k}(M)$. We can now state

\begin{thm} Let $M$ be a closed smooth Riemannian manifold of dimension $n$. Let $k$ be a positive integer, and let $2k > n$. Let $P: H^{2k}(M) \rightarrow L^2(M)$ be a $2k$th order formally self-adjoint elliptic differential operator on $M$. Let $x \in M$. Then there exists a unique weak solution of the equation $Pu=\delta_x$. 
\end{thm}

\begin{proof} First note that due to a above lemma, we have that $\delta_x \in H^{-k}(M)$. Recalling that $B(\cdot,\cdot)$ is a inner-product on $H^k(M)$ we can now invoke the Riesz Representation Theorem to conclude that there exists a function $u$ such that $B(u,v)=<\delta_x,v>$ for all $v \in H^k(M)$. Thus a weak solution of $Pu=\delta_x$ exists. 

We will now see that a weak solution is unique, Suppose that there exists a $H^k(M)$ function $u$ such that $B(u,v)=<\delta_x,v>$ for all $v \in H^k(M)$ and a $H^{k}(M)$ function $w$ such that $B(w,v)=<\delta_x,v>$ for all $v \in H^k(M)$. It would follow that $B(u-w,v)=0$ for all $w \in H^k(M)$. Setting $v=u-w$ , we have that $\beta ||u-w||^2  \leq B(u-w,u-w) = 0$, for some positive real number $\beta$. The theorem follows.
 \end{proof}

We can now state

\begin{cor} Let $(M,g)$ be a closed smooth Riemannian manifold of dimension $n$. Let $k$ be a positive integer and suppose $2k > n$. Let $P: H^{2k}(M) \rightarrow L^2(M)$ be a $2k$th order elliptic differential operator on $M$ Let $x \in M$. Let $G_x$ be a weak solution of the equation $Pu=\delta_x$. Then $G_x$ is continuous and $G_x(x)>0$.
\end{cor}

\begin{proof} Since $G_x$ is a weak solution of $P G_x = \delta_x$, we have that 
$$
\int_M P^{1/2} G_x P^{1/2} v \, dv_g = <\delta_x, v>
$$ 
for all $v \in H^k(M)$, It follows that $G_x \in H^k$, and hence is continuous. 

Now setting $v = G_x$, we have $\int_M P^{1/2}G_x P^{1/2} G_x \, dv_g = <\delta_x, G_x>$. It follows that $G_x(x)>0$. The corollary follows.
\end{proof}

\subsection{Averages of Green's functions associated with positive Green's operators are positive on symmetric sets}

\begin{thm}
Let $(M,g)$ be a closed, smooth Riemannian manifold if dimension one or higher. Let $\Omega$ be a measurable subset of $M$. Let $\chi_\Omega$ be the characteristic function on $M$ that is supported on $\Omega$.  Let $k$ be a positive integer and suppose that $P: C^\infty(M) \rightarrow C^\infty(M)$ is a $2k$th order positive elliptic differential operator on $M$. Let $G(x,y)$ be the Green's function for $P$. Then we have that  $\int \int_{\Omega \times \Omega} G(x,y) dx dy \geq 0$.
\end{thm}

\begin{proof} 
First note that since $P$ is a positive elliptic differential operator on $M$, we have $0  \leq \int_M u(x)Pu(x) dx$ for all $u \in H^{2k}(M)$.  This, in turn , allows us to write $0 \leq \int_M f(x) (P^{-1} f(x)) dx$, for all $f \in L^2(M)$. It follows then that for all $f$ in $L^2(M))$, 
$$
0 \leq \int_M f(y) (\int_M G(x,y) f(x) dx) dy.
$$
This implies that $0 \leq  \int \int_{M \times M}( f(x) G(x,y) f(y))  dx dy$ for all $f \in L^\infty(M)$, because $G(x,y)$ is integrable on $M\times M$. We can now write that 
$$
0 \leq \int \int_{M \times M} (\chi_{\Omega}(x) G(x,y) \chi_{\Omega}(y)) dx dy.
$$
This is equivalent to 
$$
0 \leq \int \int_{M \times M} \chi_{\Omega \times \Omega}(x,y) G(x,y) dx dy,
$$
where $\chi_{\Omega \times \Omega}$ is the characteristic function supported on $\Omega \times \Omega$. This in turn lets us conclude that
$$
0 \leq \int \int_{\Omega \times \Omega} G(x,y) dx dy.
$$
The theorem follows.
\end{proof}

\begin{cor} Let $(M,g)$ be a smooth, closed Riemannian manifold of dimension $n$, where $n \geq 1$. Let $P: C^\infty(M) \rightarrow C^\infty(M)$ be a positive elliptic differential operator on $M$. Let $G$ be the Green's function associated with $P$. Then we have that $\int \int_{M \times M} G(x,y) dx dy \geq  0$. Furthermore we have that $u$ is positive somewhere on $M$ if $u = P^{-1}f$ for some positive function $f$.
\end{cor}

\section{Bibliography}

\noindent [Raske 2026] D. Raske. Lowest eigenvalues and formally self-adjoint fourth order elliptic differential operators. arXiv: 2601.11882.
\vskip 0.1in
\noindent [Seeley 1967] R. Seeley. Complex powers of an elliptic operator. Proc. Symp, Pure Math. 10 (1967) 288-307.

\vskip 0.1in


\noindent [Yang 2009] P. Yang. private communication. 2009.

\end{document}